\documentstyle[amscd,amssymb,12pt]{amsart}
\input xypic \xyoption{curve}
\input epsf

\def\hcorrection#1{\advance\hoffset by #1 }
\def\vcorrection#1{\advance\voffset by #1 }

\vcorrection{-.7in}
\hcorrection{-.5in}

\newcommand{\B}[1]{{\bold#1}} 
\newcommand{\C}[1]{{\mathcal#1}} 

\theoremstyle{plain}
\newtheorem{th}{Theorem}[section]

\newtheorem{lem}{Lemma}[section]
\newtheorem{prop}{Proposition}[section]

\theoremstyle{definition}
\newtheorem{defin}{Definition}[section]

\theoremstyle{definition}

\theoremstyle{remark}
\newtheorem{rem}{Remark}[section]

\numberwithin{equation}{section}

\begin{document}

\setcounter{tocdepth}{3} 

\pagestyle{plain}
\addtolength{\footskip}{.3in}

\title[Feynman Processes]{From operads and PROPs to Feynman Processes}
\author{Lucian M. Ionescu}
\address{Department of Mathematics, Illinois State University, IL 61790-4520
}
\email{lmiones@@ilstu.edu}
%

\begin{abstract}
Operads and PROPs are presented, together with examples and applications
to quantum physics suggesting the structure of Feynman categories/PROPs and 
the corresponding algebras.
\end{abstract}

\maketitle
\tableofcontents


\section{Introduction}\label{S:intro}
The concepts of operads and PROPs allow to formalize certain types of
algebraic structures, starting with the classical 
commutative, associative, and Lie algebras, 
but most notably those appearing in the 
quantum mathematical-physics.

The more general framework is that of multi-categories \cite{Lesner},
which are higher-dimensional categories with morphisms 
(``multi-arrows'', ``coupons'' etc.).
Composing such ``multi-arrows'' (with several inputs) needs
prescriptions for composing ``composable arrows''.
The asymmetry n-inputs/1-output is not beneficial,
allowing for the specialization of MacLane's concept of 
PROP (Product Operations and Permutations) to
P. May's combinatorial definition of operads,
which hides the true meaning, 
being quite obscure at first sight \cite{DF}.

We will follow the short introduction \cite{DF}, 
including comments referring to their applications 
to (quantum) mathematical-physics.

Starting from PROPs, 
which are monoidal categories enriched over another category
(Hom spaces are not just sets, being endowed with additional structure),
operads will appear naturally as ``half-PROPs''.
For pedagogical purposes we will introduce first the version 
without the equivariant action of the permutation groups 
(``non-$\Sigma$'' operad/PROP - see \cite{MSS}); 
this corresponds, for example at the graphical level,
to considering unlabeled graphs.
Various examples are provided, 
and the final section is concerned with the various quantum models viewed as 
algebras over a corresponding PROP \cite{V}.

\section{PROPs}
The intrinsic definitions will be followed by explicite descriptions.
\begin{defin}
A PRO in $C$ over $\C{N}$ is a strict monoidal category with objects of $\C{N}$
and Hom spaces objects of $\C{C}$.
\end{defin}
We will be concerned with the case when $\C{N}=(N,+)$,
the monoid of natural numbers,
i.e. the strict monoidal category is generated by one object.

\subsection{PROs}
Explicitly, for each pair of natural numbers $(n,m)$, 
an object $P(n,m)$ in $C$ is given, 
representing the corresponding Hom space,
together with two types of composition \cite{M}, p.37:
$$\otimes:P(m_1,n_1)\otimes ...\otimes P(m_s,n_s)
\to P(\sum m_i, \sum n_i),$$
referred to as ``horizontal composition'', and
$$\circ: P(m,n)\otimes P(n,k)\to P(m,k),$$
called ``vertical composition''
\footnote{Convention: ``time'' vertical and ``space'' horizontal}. 
The two compositions are required to be associative and unital.
The monoidal unit $1\in P(0,0)$ should not be confused with 
the identity morphism $I\in P(1,1)$.

In more detail, an alternative explicite definition of a PROP is taken from \cite{DF}.
\begin{defin}
The data for a $\C{C}$-PROP over $N$,
i.e. with $N$ as a category of objects and $\C{C}$ as a category of morphisms
consists in:

1) A map $P:N\times N\to \C{C}$;

2) Composition maps: $\circ:P(m,n)\otimes^{op} P(n,k)\to P(m,k)$;

3) Tensor product maps: $\otimes:P(m,n)\otimes P(k,l)\to P(m+k,n+l)$;

4) A composition identity $j\in P(n,n)$ and possibly a tensor unit $1\in P(0,0)=k$,

satisfying the natural associativity and compatibility conditions.

Braided or symmetric PROPs, also called $\Sigma$-PROPs \cite{MSS},
have in addition braidings 
$\sigma_{m,n}\in P(m+n,n+m)$, inducing a structure of a
braided/symmetric tensor category.
\end{defin}

The compatibility between the two operations,
similar to the one in a monoidal category
$$(f\otimes g)\circ(f'\otimes g')=(f\circ f')\otimes (g\circ g')$$
requires that the category $\C{C}$ be symmetric (see \cite{DF}, p.6).
\begin{rem}\label{R:functor}
It is useful to think of a category as a ``functor'' \cite{DF}, p.2:
$$Hom:Ob(\C{N})^{op}\times Ob(\C{N})\to Sets,$$
and visualize it as a multi-graph (multiple arrows between labeled vertices).
Then an enriched category is just the general case:
$$Hom:Ob(\C{N})^{op}\times Ob(\C{N})\to Ob(\C{C}).$$
In view of the applications to physics,
such a structure encodes the possible {\em states and transitions} (paths)
of a dynamical system, consisting of interacting subsystems.

Being concerned with FPI method,
a {\em Feynman Category} will be defined as a PROP of dg-coalgebras,
for reasons to be explained later on.
\end{rem}

\subsection{PROPs}
A PROP is an equivariant version of a PRO,
i.e. the Hom spaces, objects of $\C{C}$, are G-spaces
under the symmetric action of the 
corresponding permutation groups: $P(m,n)\in Ob(\C{C})$ 
is endowed with a left $S_m$ and right $S_m$ action,
compatible with the above composition operations.

\subsection{The basic example}
The basic example is {\em endomorphism PRO} 
associated to a k-vector space $V$:
$$End(V)(n,m)=Hom_k(V^{\otimes n}, V^{\otimes m}),$$
together with the tensor product $\underset{k}{\otimes}$ and
composition of k-multilinear mappings.

The following remark provides also examples of PROPs.
\begin{rem}
To understand the above ``natural'' conditions,
note that a tensor category $\C{C}$ is tautologically a 
$\C{C}$-PROP over $\C{C}$.
Therefore the full tensor sub-category generated by an object $``1''$
is a $\C{C}$-PROP over $N$.

One can slightly generalize the concept of PROP by taking
$Z$ as a model for objects; for instance $\C{C}=SU_2(C)-mod$ (spin)
would be such an example, were not for the more complicated
fusion rule for spins.
\end{rem}
Operations with PROPs will be considered after specializing PROPs to operads,
in order to consider some familiar examples first.

\section{Operads}
``An operad is simply the $n-to-1$ part of a PROP ...'', \cite{DF}, p.8:
$$O(n)=P(n,1).$$
Note that an operad is not a genuine subobject,
rather it is a ``truncation/restriction'' of a PROP.

\subsection{Restricting a PROP}
\begin{defin}
The collection of morphisms $O(n)=P(n,1), n\in N$ of a 
$\C{C}$-PROP over $N$, is called an {\em operad}.

If the PROP is a symmetric/braided PROP, then so is the operad.
\end{defin}
We will give an explicit definition bellow.
After ``hiding'' the 2nd argument of $Hom$,
the (P. May's) definition becomes quite ``encrypted'' 
by not distinguishing the two compositions $\circ$ and $\otimes$.
%
\begin{defin}
An {\em operad} $\C{O}$ in $\C{C}$ (or of $\C{C}$s) 
consists of a collection of objects of $\C{C}$
$$\C{O}(n)=P(n,1)\in Ob(\C{C}), \quad n\in N$$
together with composition operations
$$\circ:O(n)\times (O(k_1)\otimes ...\otimes O(k_n))\to O(k_1+...+k_n),$$
satisfying the (natural) {\em associativity} and {\em unital} conditions.
\end{defin}
Denoting by $\theta_n$ a generic operation of $O(n)$,
the associativity condition is:
$$\theta_n\circ (\theta_1\circ
(\theta_{1,1}...\theta_{1,k_1}),...,\theta_n\circ(\theta_{1,1}...\theta_{1,k_1}))
=(\theta_n\circ (\theta_1...\theta_n))\circ
(\theta_{1,k_1}...\theta_{n,k_n}),$$
and the {\em unital} requirement is: 
there is an element $j\in \C{O}(1)$, called unit, such that
$$\theta_n \circ(j,...,j)=\theta=j\circ \theta_n, \quad \theta_n\in O(n).$$
In conclusion, being ``economical'' has a too high price
we should not pay, and avoid the specialization.

\subsection{The basic example}
The basic example is the endomorphism operad $End(V)(n)=Hom(V^{\otimes n},V)$,
together with composition of k-multilinear mappings;
properly speaking it is an algebra over the tree operad,
as it will be explained later on.

If we ``add'' the action of the symmetric group to operads of sets,
and require that composition of operations of $\C{C}$ 
are equivariant morphisms, we obtain the notion of {\em symmetric operad}.
The endomorphism operad is a symmetric operad.

\subsection{PROP generated by an operad}
Given an operad $O$, we can define a PROP $P$ having $O$ as its
restriction $O(n)=P(n,1)$ \cite{DF}, p.8.
In order to define $P(m,n)$,
recall that for any $l$-partition of $m=m_1+...+m_l$,
the tensor product is a map:
$$\otimes: P(m_1,1)\otimes ...\otimes P(m_l,1)\to P(m,l).$$
Assuming $\C{C}$ has direct sums, and since $P(m_i,1)=O(m_i)$, set:
$$P(m,l)=\oplus_{(m_1,...,m_l)|\sum m_i=m}
O(m_1)\otimes ... \otimes O(m_l),$$
together with tensor products defined by the canonical inclusions and 
compositions defined component-wise (see \cite{DF}).

The above PROP is {\em  the PROP generated by the operad $O$}
and denoted by $P=\C{P}(O)$.
Together with the {\em underlying operad} of a PROP is $O=\C{U}(P)$,
it has the following properties:
$$\C{U}(\C{P}(O))=O, \quad \C{P}(\C{U}(P))\sim ?.$$
The standard example is again the endomorphism PROP/operad;
in this case $\C{P}(Hom(V^\bullet,V))\cong End(V^\bullet)$.

\section{Representations of PROPs and operads}
This are linearized versions of morphisms of PROPs and operads.

\subsection{Morphisms of PROPs}
A morphism of PROPs is a (strict) monoidal ``functor'' (see \ref{R:functor})
$\rho:P_1\to P_2$ between two such enriched monoidal categories:
$$\rho_{Ob}:\C{N}_1\to \C{N}_2,$$
$$\rho_{Hom}:P_1(A,B)\to P_2(\rho_{Ob}(A),\rho_{Ob}(B)), A,B\in \C{N},$$
intertwining the additional structures of the two $Hom$ spaces.

Being strict means that $\rho(A\otimes B)=\rho(A)\otimes \rho(B)$
and $\rho(\theta_1\otimes \theta_2)=\rho(\theta_1)\otimes \rho(\theta_2)$,
where the subscripts were omitted for convenience. 

As an example take $P_2=Vect_k$, the PROP associated to the
strict tensor category $Vect_k$,
and $P_1=P$ a $\C{C}$-PROP over $N$.
Then a morphism $\rho:P\to Vect_k$ is determined by 
choosing a vector space $V=\rho_{Ob}(1)$,
determining therefore $\rho_{Ob}(n)=V^{\otimes n}$,
together with a family of morphisms 
$$\rho_{m,n}:P(m,n)\to Hom(V^{\otimes m},V^{\otimes n})$$
which is ``functorial'', i.e. compatible with horizontal composition $\circ$,
and ``monoidal'', i. e. compatible with the tensor product.

Since the range of the morphism is the tensor full subcategory
$(V^\bullet, End_k(V^\bullet)$,
isomorphic to the endomorphism PROP $\C{E}nd(V)$,
the above morphism of PROPs co-restricts to the 
following morphism of PROPs:
$$\rho:P\to \C{E}nd(V).$$

\subsection{Representations: algebras over a PROP}
In the above example, 
elements of $P(m,n)$ act as multilinear operations on $V$.
\begin{defin}
A vector space $V$ together with a morphism $\rho:P\to \C{E}nd(V)$
is called a $P$-algebra.
\end{defin}
For a concrete example, 
define $P(1,1)=G$ a group with all other $Hom$ spaces empty,
to view a group representation $\rho:G\to End(V)$ as a morphism of PROPs
$\rho:P\to \C{E}nd(V)$.
In general, a $P$-algebra will be also called a representation
of the PROP $P$.

Assuming that both PROPs have the same objects $\C{N}_1=\C{N}_2=:\C{N}$,
and that the morphism is ``over the identity map'',
i.e. $\rho_{Ob}=Id_{Ob(\C{N})}$, 
then such a morphism can be described as follows.
\begin{defin}
A morphism of PROPs $\rho:P_1\to P_2$ over $N$ is 
a family of mappings $\rho_{A,B}:P_1(A,B)\to P_2(A,B)$ 
compatible with the two operations $\circ$ and $\otimes$.
\end{defin}
If $\C{N}\cong (N,+)$ is generated by one element $1$,
then, by the above arguments, the morphism of PROPs
is equivalent to the $\C{C}$-representation
$$\rho:P\to \C{E}nd_\C{C}(A),$$
where $A=\rho(1)$; it will be called
a $P$-algebra in $\C{C}$.

A (1+1)-TQFT, to be defined later, 
provides such an example.

\subsection{Morphisms of operads}
If we restrict to operads, an $O$-algebra 
is a family of n-ary morphisms:
$$\rho_n:O(n)\to Hom(V^n,V),$$
compatible with composition and tensor product.
Since elements of $O(n)$ act as n-ary (k-linear0 operations on $V$,
it is natural to call such a structure 
a representation of the operad $O$.
\begin{rem}
The categorical point of view may be traded for a classical terminology
by formally adding all the morphisms of an operad (or PROP) into
one graded (bigraded) ring.
As an important example, we point to Hopf algebras of Feynman graphs,
and to graph complexes in general.
In these cases the classical terminology applies:
subobjects, ideals and quotients etc.,
in spite of the prefered categorical description.
\end{rem}

%
\section{Operations with PROPs and operads}
To simplify matters, we will work with k-linear PROPs and operads,
where $k$ is a field.

\subsection{Free operads}
A free operad solves a universality problem, as expected:
given a collection of sets $X=\{X(n)\}$ find an operad of vector spaces
$\C{F}(X)$ and an inclusion map $j:X\to \C{F}(X)$, 
such that for any operad of vector spaces $O$ and mappings of sets
$\rho_n:X(n)\to O(n)$ there is a unique lifting of $\rho$ 
to a morphism of operads $\C{F}(X)\to O$,
i.e. the pair $(\C{F}, j)$ is an initial object in the appropriate 
category of operadic valued mappings.

The above problem has a solution called the {\em free operad}
generated by the sets $X(n)$.
Applying the result to $O=\C{E}nd(V)$, 
the endomorphism operad, we conclude that any choice
of a set $X(n)$ of n-ary multiplications on $V$ defines 
an $\C{F}(X)$-algebra structure on $V$.

The existence of the free operad of vector spaces can be understood
as a consequence of the existence of a free operad of sets generated
by the collection $X$, together with the existence of a (free)
vector space with a given basis.

\subsubsection{Trees and forests}
The free operad generated by a given set,
is the operad of the free PROP of forests, colored by the elements
of the given sets $X(n)$.

To give the additional details, at this point we adopt a presentation
which avoids the embedding of trees in a plane (compare with \cite{DF}, p.11),
by considering {\em labeled trees} and {\em forests}.
These may be realized, for instance, 
as sequences of order preserving maps of finite sets.

A {\em planar tree} $t\in T(n)=F(n,1)$ is a tree together with an order of 
its terminals (n-points boundary vertices), 
usually given by a labeling of its terminals.
A {\em k-forest} $t_1 ... t_k \in F(\sum_{i=1}^k n_i,k)$,
is an ordered concatenation of trees, 
together with an order on the boundary vertices,
possibly different from the concatenation of the orders of the
constituting trees.
The roots form the {\em outgoing boundary} $\partial^+F$,
while the terminals form the {\em input boundary} $\partial^-F$.

The tensor product is concatenation, 
with the natural induced order on the boundary,
while composition of two compatible forests $F\circ G$ 
is the gluing operation (identification) of the outgoing boundary
of one forest with the input boundary (terminals) of the other
forest, uniquely determined by the order on the respective boundaries.

Let $F(m,n)$ denote the set of all forests with $m$ inputs and $n$ outputs.
\begin{lem}
$F$ is a symmetric PROP of sets over $N$, 
with identity the unique tree $I\in F(1,1)$
and braiding represented by the transposition $(12)$ viewed
as a forest in $F(2,2)$.
\end{lem}
Note that the forest $(12)$ is different from $I\otimes I$.

Forgetting the order on the boundary vertices by identifying
the corresponding forests, leaves us with the underlying PRO $\bar{F}$.

Its underlying operad is the {\em operad of trees}.
It is generated as an operad of sets by the n-to-1 corollas $X(n)=\{(n,1)\}$,
i.e. it is the {\em free PROP of sets} generated by the collection of sets $X$
(graded set).
\begin{prop}
The PRO $\bar{F}$ is the PRO generated by the free operad of trees
and also the free PRO generated by the corollas.
\end{prop}
Including the order on the vertices of the corollas
$X(n)=\Sigma_n$ one obtains the corresponding free symmetric operad.
To obtain the free symmetric PROP one must apply also the shuffles
to concatenations of ordered trees.

\subsubsection{Colored forests}\label{S:colors}
Now consider the sets $X(n)$ as colors for the n+1-valent internal vertices
of the trees of the forests of $F$. 
The {\em X-colored forests} form a PROP $F(X)$: it is the free PROP of sets
generated by the sets $X$.

As a special case, $F$ may be thought as forests generated by 
corollas of just one color.

If $X(n)=\emptyset$ for some $n$, then these $n+1$-valent nodes are
forbidden.

We can also linearize and consider the corresponding PROPs of vector spaces,
generated by $F(X)$ over $k$, and denoted $F(X;k)$.
Now a morphism of PROPs of vector spaces $\rho:F(X;k)\to P$,
is determined by the a family of maps of sets $\rho_n:X(n)\to P(n,1)$.
This shows that the underlying operad of $F(X)$ is 
the free operad we were looking for;
its elements are linear combinations of trees with 
ordered boundary, colored by elements of $X$.

\begin{rem}
Another important example to be discussed later on,
is the PROP of forests decorated by 1-PI diagrams,
isomorphic to the PROP of Feynman diagrams.
\end{rem}

\subsection{Ideals and quotients}
Let $P$ be a PROP of vector spaces.
An ideal of $P$ is a subcategory of the enriched tensor category
having the same objects as $P$, 
i. e. a subcategory which is still a PROP.
Explicitly, we have the following
\begin{defin}
An {\em ideal} of a $\C{C}$-PROP $P$ over $N$ 
is a collection of normal subobjects $I(m,n)$ of $P(m,n)$ closed
under composition and tensor products.
\end{defin}
One may consider in what follows PROPs of vector spaces for simplicity.
\begin{rem}
The term {\em sub-PROP} would be more appropriate; 
an ideal should rather be closed under products with 
general elements of the PROP.
\end{rem}
It is immediate that the family $(P/I)(m,n)=P(m,n)/I(m,n)$
is again a PROP.

If $I\subset P(O)$ is an ideal of the PROP generated by the operad $O$,
then the underlying operad of the quotient PROP is called the
quotient PROP of $O$ by $I$.
Note that we cannot speak of the ideal of an operad directly,
but we can consider a suboperad $S\subset O$, 
i. e. a family $S(n)\subset O(n)$ of normal subobjects closed under composition.
Then $P(S)$ is an ideal of the PROP $P(O)$ and
the quotient $P(O)/P(S)$ is the quotient operad.

Ideals in PROPs may be obtained as generated from a collection of subsets
of a concrete PROP, as the smallest ideals containing the given sets.

\subsection{Duality and cooperads}
There is a natural notion of a dual PROP $P^{Op}$,
modeled after the dual of a category.
The underlying operad of the dual PROP is called a {\em co-operad}.
The corresponding concepts can be readily obtained by duality.

We will stop with generalities and focus on some basic examples.

\section{Classical examples of operads}
Again we will follow closely \cite{DF}.

\subsection{The operad {\bf Assoc}}
Let $X(2)=\{ (1,(12)) \}$ and $X(n)=\emptyset$ for $n\ne 2$.
The (non-sigma) operad $F(X)$ is the operad of binary trees
(trivalent planar trees); its representations are algebras with 
one binary operation.

The ideal of the free operad, with 1 generator:
$$\C{I}=<\circ>,\quad \circ=\circ_1-\circ_2.$$
encodes the associativity constraint.
The corresponding quotient operad $\C{F}/\C{I}$ will be denoted by ${\bf Assoc}$.
The {\bf Assoc}-algebras are the associative algebras.

Alternatively, one can view associative algebras together with
the symmetric action as algebras over the symmetric operad
corresponding to the free operad generated by the labeled 2-corolla:
$X_\sigma(2)=\{ (12), (21)\}$ and ideal generated
by $(12)3-1(23)$.

\subsection{The operad {\bf Comm}}
It is also a quotient of the associative operad ${\em Assoc}$
by the relation $12-21$ (it is not a planar operad).
The {\bf Comm algebras} are the commutative and associative algebras.

\subsection{The operad {\bf Lie}}
It is the quotient of the free operad $F(X)$ by the 
anti-commutative condition $(12)+(21)$
\footnote{Different symmetric character used to define the braiding.}
and the Jacobi identity (``infinitesimal associativity'')
$(12)3-1(23)+2(13)$.
The {\bf Lie}-algebras are, of course, the Lie algebras.

\subsection{The operad {\bf Poisson}}
The operad Poisson involves two operations $X(2)=\{m,[,]\}$,
which may be interpreted as colors of the 2-corolla (see \ref{S:colors}).
The other sets $X(n)$ are empty, as before.
The corresponding free operad $F(X)$ is operad of binary trees
colored by the collection $C$ of colors.

The {\bf Poisson} operad is the quotient of $F(X)$
by the ideal generated by the following constraints:
associativity and commutativity for $m$,
anti-commutativity and Jacobi identity for $[,]$,
together with the compatibility between the two:
$$m([1,2],3)-[1,m(2,3)]+m(2,[1,3]),$$
or, written as a relation, 
with the usual notation $a=1, b=2, c=3$ and $m$ concatenation:
$$[a,bc]=[a,b]c+b[a,c],$$
i.e. $ad_a=[a,\cdot]$ is a derivation.

\section{Examples of PROPs}
It was briefly mentioned that Feynman diagrams form a PROP,
or alternatively, as a global classical algebraic structure,
a Hopf algebra \cite{CK} 
(more precisely dg-Hopf \cite{I-pqft,I-cfg}).

\subsection{Feynman PROPs and Feynman categories}\label{S:FC}
Indeed, Feynman diagrams with external legs partitioned into
two sets, input and output, can be viewed as elements of a PROP;
e.g. the tad-pole $tp\in FD(1,2)$ etc.
They are in a 1:1 correspondence with the forests colored by
connected 1-PI (particle irreducible) Feynman diagrams, 
and therefore form a PROP as explained above in \S\ref{S:colors}).
The details can be found in \cite{CK},
but only at the level of a PROP of vector spaces.
The coproduct structure and the resulting canonical antipode
represents an additional structure on the vector spaces ($\Sigma$-modules)
$FD(m,n)$; together with the compatible homology differential \cite{I-cfg}
it is a PROP of dg-coalgebras.
\begin{defin}
A PROP of dg-coalgebras is called a {\em Feynman Category}.
If the PROP is over $N$, it is also called a {\em Feynman PROP}.
An algebra over a Feynman Category is called a {\em Feynman Process}.
\end{defin}
Various examples can be obtained by considering the so called
$\phi^3, \phi^4, ...$ theories. 
These are classes of graphs
obtained by expanding Gaussian integrals using Wick's Theorem.
More generally, the Wick's Matrix Theorem can be used to obtain
classes of Feynman diagrams (Feynman Categories) from
{\em Matrix Models} (\cite{F-MM} etc.).
There are important correspondences between such classes of 
Feynman diagrams and other PROPs,
notably involving {\em Riemann Surfaces} \cite{F-MM}.
Additional details will be provided in \S\ref{S:appl}.

\subsection{PROPs and ``bi-operads''}
Intuitively speaking,
PROPs contain multiplications and comultiplications,
so it is natural to compare them with bialgebras,
built of an operad compatible with a co-operad.

Given two operads $A$ and $B$, one can build a PROP $P=B^*A$,
as the smallest PROP containing $F(A)$ and $F(B)^*$.
It is expected to be freely generated by the element obtained 
by requiring that all ``annihilation operators'', i.e. multiplications
must precede creation operators (coproducts).
Related topics can be found in \cite{M}.

The graphs defined by Kontsevich in \cite{K-formal} form a Feynman PROP $K$,
by rephrasing the dg-coalgebra structure from \cite{I-cfg}.
We conjecture to be the free PROP generated by the forest PROP
as sketched above: $K=F^*F$.
The meaning within the framework of quantization/renormalization
(OPEs etc.) will be investigated elsewhere.

\section{Higher operads: homotopy algebras}
Various structures involving several operations and compatibility conditions
can be unified as an algebra associative modulo higher brackets.

A {\em homotopy Lie algebra} is a graded vector space $H$ together 
with n-ary operations $m_n=[\ ,...,\ ]:V^{\otimes n}\to V$, $n=1,2,...$,
satisfying the graded equation $Q^2=0$ where 
$$Q=[\ ]+[\ ,\ ]+...$$
and the operadic composition (or rather the composition of the generated PROP)
is understood (Gerstenhaber composition in classical terms).

With the usual notation $d=[\ ]$, the first few equations are:
$$d^2=0, d[\ ,\ ]=[d\ ,\ ]+[\ ,d\ ], ...$$
The next equation, in degree 2, is Jacobi identity (Lie associativity) 
``modulo homotopy'' represented by the next bracket $m_3$.
As a basic example, a DGLA is a homotopy algebra for which 
$m_1=d, m_2=[\ ,\ ]$ while all the other operations are trivial.

Homotopy algebras are algebras over a quotient of the tree operad, 
viewed as a symmetric operad with anti-symmetric braiding; 
its tensor algebra is the algebra over the forest PROP.
One can check that the composition is the one of the 
standard PROP generated by a given operad. 

The principal types of homotopy algebras are, as in the classical case,
associative, commutative and Lie \cite{K-formal,GK}.
The Lie and commutative operads are dual to one another, 
while the associative operad is self-dual.

For more details see also \cite{KSV}.

\section{Applications: TQFTs, CFT, ST and homotopy algebras}\label{S:appl}
The more sophisticated examples of representations of Feynman Categories (FC), 
modeling quantum interactions/transitions, 
are TQFTs, CFT, ST and involve a generalization of the usual classical 
finite operads/PROP, i.e. $O(n)\ne 0$ for possibly all $n$: homotopical algebras.
We will follow \cite{V} regarding the main examples.

What we call FC corresponds to Voronov's ``state space'' of the 
{\em geometric field theory}, which has the flavor of a 
``geometric category'', at least at the level of its graphical representation,
if not because it represents the {\em external causal structure}.
It stands for a more general model then the classical space-time,
taking into account economically the {\em actually existing} subsystems, 
rather then the potential ones, and leading to fixed empty space, etc..

We call the representations of a Feynman category, 
i.e its algebras over the corresponding PROP, 
a Feynman Process, to point towards the state-transition
interpretation ``a la Markov''.

\begin{rem}
We will focus on the 2-real dimensions case of the geometric causal
structure, since what is being represented is the generic ``parallel
versus sequential'' relationship/interaction, 
which classically appears as space versus time, 
including the complex 1-dimensional case of Riemann surfaces. 
Ultimately, it is the functorial representation of how the complex information
is processed (computer science perspective) and interacts (quantum physics perspective), 
as a 1D-complex flow instead of the classical 1D-real time flow. 
The internal DOFs appear as part of an algebraic structure ``on top'' of
the geometric Feynman category, in the spirit of Feynman rules.
\end{rem}
In this vein, it is not surprising that the algebraic string
background of a string theory involves the Virasoro semi-group of cylinders
acting on some Hilbert space of internal states.

\subsection{TQFTs}
TQFTs are representations of cobordism categories \cite{A}.

\subsubsection{(1+1)-TQFTs}
A (1+1)-TQFT is a functor having as a geometric category 
the Feynman Category of surfaces with boundaries a disjoint union of circles, 
viewed as cobordisms under a suitable equivalence relation
ensuring a well defined gluing as composition. 
As in \cite{V}, 
we will call this the {\em topological Segal PROP} {\em Top-RS}, 
since we forget the complex structure of RS.

It is a tensor functor valued in $\B{C}$-Vect (turns out to
be k-Hilbert: vector spaces with self-duality):
$$Top-RS(m,n)\ni \Sigma \mapsto V\in \B{C}-vect,$$
satisfying natural axioms: 1) topological invariance, 2) permutation
equivariance, 3) factorization property, 4) normalization; these are
condition ensuring that the algebraic structure of k-Vect
``matches'' the topological equivalence and symmetry structures 
(see \cite{V} for additional details and pictures).
\begin{prop}
An algebra over the topological Segal PROP is a (1+1)-TQFT. 
\end{prop}
\begin{pf}
The statement follows from an inspection of Atiyah's 
definition of a TQFT \cite{A}.
\end{pf}
A folklore result relates the set-theoretic description and the
categorical description of the same physical theory \cite{Kock}. 
\begin{th}
A (1+1)-TQFT is equivalent to a Frobenius algebra.
\end{th}
Recall that a Frobenius algebra is a commutative algebra with unity and an invariant
non-degenerate symmetric bilinear form invariant with respect to
multiplication:
$$<ab,c>=<a,bc>$$
which has an adjoint (copairing) $C\to V\otimes V$.

This is the associativity of ``string interaction''. 
The copairing plays the role of the existence of a basis in the state space
$V=F(O)$.
This case is a typical example when the abstract, non-linear algebraic object gets represented
in the realm of linear spaces \cite{I-note}. 

An important subobject is the {\em TQFT at the tree level}, i.e. the
restriction of the functor to the underlying operad
(genus zero Top-RS).
\begin{th}\label{T:tree-level1}
The underlying operad of Top-RS (tree level), 
consists in the genus zero topological Riemann surfaces.

A (1+1)-TQFT at the tree level is equivalent to a commutative
algebra with unity.
\end{th}
This is in a way expected, since at tree level one should recover
the classical limit, i.e. the observable algebra $A$ 
without quantum corrections due to ``interference'' and loops,
and since a commutative algebra can be represented as the observables on a compact space
$\C{F}(Spec A)$. 
being topological, the dynamics is trivial and the Poisson bracket is trivial.

Comparing the two cases, to extend the functor from the operad (tree level) to the PROP
amounts to introducing the invariant pairing.

\subsubsection{(0+1)-TQFTs}
The cobordims with points as boundaries are built out of 
compact intervals and circles.
This leads to categories of tangles and their representations
\cite{Turaev}.

\subsection{Conformal Field Theory}
A conformal field theory (CFT) is a functor similar to a 2D-TQFT:
$$CFT: \sigma \mapsto |\Sigma>,$$
except:

1) The Segal PROP consists of Riemann surfaces bounding holomorphic disks;

2) The functor $|\Sigma>$ depends smoothly on the RS;

3) Topological invariance is replaced by {\bf conformal invariance};

4) The topological gluing is now supplemented by gluing the complex structures;

5) Normalization is now:
$$O\mapsto id_H,$$
i.e. a cylinder of zero length is mapped to the identity map of the
Hilbert space.

In other words, a CFT is a smooth mapping:
$$\C{P}_{m,n}\to Hom(H^m,H^n).$$
where $\C{P}_{m,n}$ is the moduli space of RS, bounding m negatively
oriented disks, and n positively oriented disks.
%

\subsection{String Theory and Homotopy Lie Algebras}
String theory is usually associated with its sigma-model incarnation,
where the action is obtained by embedding Riemann surfaces in an ambient
space (manifold). 
To avoid a pre-assumed space-time or target Riemannian manifold,
an algebraic formulation proceeds as follows.

\subsubsection{String backgrounds} 
A {\em string background} (SB) is a correspondence:
$$SB: C_\bullet\C{P}_{m,n}\to Hom(H^m,H^n)$$
where now $(H,Q)$ is a graded complex ($Q^2=0$). 
In other words, an SB is a morphism of PROPs of complexes over the Segal category, 
i.e. between the Segal DG-PROP and the Hochschild algebra over a forest DG-PROP 
associated to the complex $(H,Q)$.
\begin{rem}
Note this is a half-step towards the structure of a Feynman Category
which is a PROP of dg-coalgebras (DGLAs).
\end{rem}
Here $C_\bullet\C{P}_{m,n}$ is 
the complex of chains of the moduli space of RS \cite{KSV}. 

As a consequence, a SB must satisfy the following axioms:

1) Smoothness: the mapping is smooth 
(the moduli space is a finite dimensional manifold);

2) Equivalence: the mapping is equivariant with respect to permutations; 

3) Factorization: the mapping respects sewing the RS corresponding
to the composition of two chains;

4) Homogeneity and $Q-\partial$-invariance: the
mapping is a chain map (morphism of complexes); 

5) Normalization: the R Sphere with two unit disks cut is mapped to the identity map.

In the case of SB, conformal invariance is not factored out, as for CFT:
the zero width cylinders together with diffeomorphisms of the circle
define the {\em Virasoro semi-group of cylinders} \cite{V}, p.7.
If $v$ is a vector field on the circle and $exp(v)$ the corresponding diffeomorphism,
then $exp(tv)$ is a 1-chain $B(tv)$ of degree -1 acting on $H$:
$$cylinder(v)\mapsto B(tv)=exptT(v):H\to H.$$

Descending to homology, one gets the associated (1+1)-TQFT.
\begin{th}
Forgetting the complex structure of the Segal PROP and taking the
homology of $H$, yields a TQFT. Thus the homology of $H$ has a
natural structure of a Frobenius algebra.
\end{th}
This is a string theory ``incarnation''
of the de Rham / Hodge representative correspondence.

\subsubsection{Higher brackets}
The genus zero RS ($M+1$-punctured Riemann spheres), 
form an operad $O(m)=\C{M}_{g=0}(m,1)$.
In a way similar to the TQFT case,
restricting a SB to the genus zero operad underlying the Segal DG-PROP
yields a homotopy Lie algebra.

The higher brackets $[\ ,...,\ ]=SB(O(m))$ are defined using 
additional data \cite{Z} leading to {\em closed string-field theory}.
\begin{th}\label{T:tree-level2}
The underlying operad of the Segal DG-PROP (tree level), 
consists in the genus zero punctured Riemann spheres.

A (1+1)-SB at the tree level is equivalent to a homotopy algebra.
\end{th}
Yet the smoothness condition is too week,
and the natural condition to be imposed on a SB is analyticity.

\subsubsection{Holomorphic SB and homotopical commutative algebras} 
The duality between the homotopy Lie and commutative operads 
corresponds at the level of string theories to passing 
from {\em left-right movers} to {\em chiral case},
roughly corresponding to passing from 
smooth ``operator-valued differential forms'' on the moduli space 
to holomorphic ones (\cite{V}, p.10; \cite{KSV}).

\section{Conclusions: Feynman Processes}
The above classical examples of mathematical models of quantum theories
show clearly a pattern: modeling quantum processes by states and transitions 
in the spirit of Markov processes, except graded by complexity (``perturbatively'')
and providing complex amplitudes of probabilities, as expected.

The heritage of the classical deterministic-mechanicist view is represented by the ambient 
space-time manifold used to define the action in sigma models.
As pointed out in several occasions, our goal is to formulate the intrinsic theory first,
and then look for representations, constructed algebraically as for string backgrounds
or using an ambient space-time as in sigma models.

The basic structure needed in such a intrinsic formulation was introduced in 
S\ref{S:FC}: Feynman categories are dg-coalgebra PROPs and their representations, 
the Feynman Processes are Feynman algebras.
The key additional feature not used at the level of Hopf algebras of Feynman graphs,
is the homology differential compatible with the coalgebra structure yielding the
scalability of the theory.
This enables the interpretation of the ``perturbative'' approach as rather 
a resolution of the quantum system, in the sense of homological algebra \cite{I-FL}
(see also the {\em Quantum Dot Resolution} \cite{I-DWT}).

This leads to another important point learned from the above examples \S\ref{S:appl}, 
that the homology of the theory corresponds to the tree level and 
classical solutions, while the theory at chain level includes the quantum corrections.
A more detailed study of this point of view, is deferred to another article,
incorporating the cohomology of Feynman graphs from \cite{I-cfg}.

More importantly still is the deep reason for the absence of a 
global ambient space-time in such an intrinsic approach via 
Feynman categories and algebras,
or the lack of a global state space (manifold) as in classical mechanics, general relativity etc.
The meaning comes from the computer science interpretation, that quantum interactions
are quantum communications, i.e. processing quantum information \cite{I-FL,I-DWT}:
a ``space-time structure'' is only a ``decomposition'' of the 
computing process in sequential and parallel operations,
which leads in the large scale at the familiar space-time manifold model,
and which manifests itself as an {\em operator product expansion}
at the level of chiral algebras.
This will be explained in detail elsewhere \cite{I-QID}.

\section{Acknowledgments}
I would also like to thank for the research and lodging conditions
while visiting I.H.E.S., where some of the present ideas matured and contributed
to a better ``big picture'' of the subject matter in the second talk.



\end{document}